\documentclass{kms-b}
\issueinfo{}
  {}
  {}
  {}
\pagespan{1}{}
       
\usepackage{graphicx}
\allowdisplaybreaks
\usepackage[utf8]{inputenc}
\theoremstyle{plain}
\usepackage{amsmath}
\usepackage{amsfonts}
\usepackage{amssymb}
\usepackage{graphicx}
\usepackage{subcaption}
\usepackage{booktabs}
\usepackage[all,cmtip]{xy}
\usepackage{algorithmic}
\usepackage[linesnumbered,ruled]{algorithm2e}
\usepackage[hmargin=1.5cm,v margin=2.5cm]{geometry}
\usepackage{vaucanson-g}
\usepackage[T1]{fontenc}
\usepackage{float}
\usepackage{tikz}
\usetikzlibrary{arrows}
\newtheorem{Theorem}{Theorem}[section]      
\newtheorem{Definition}[Theorem]{Definition}
\newtheorem{Conjecture}[Theorem]{Conjecture}
\newtheorem{Proposition}[Theorem]{Proposition}
\newtheorem{Remark}[Theorem]{Remark}
\newtheorem{Proof}[Theorem]{Proof}
\newtheorem{Example}[Theorem]{Example}
\newtheorem{Corollary}[Theorem]{Corollary}
\newtheorem{Lemma}[Theorem]{Lemma}

\begin{document}

\title[A New Hyperbola based Approach to factoring Integers]
{A New Hyperbola based Approach to factoring Integers}

\author[Bansimba]{Gilda Rech Bansimba$^{1}$}
\address{Gilda Rech Bansimba \\ Department of Mathematics and Computer Science\\ Universit\'e Marien Ngouabi \\ Facult\'e des Sciences et Techniques \\ BP: 69, Brazzaville, Congo \\
	https://orcid.org/0000-0001-9213-6930}
\email{bansimbagilda@gmail.com}

\author[Babindamana]{Regis Freguin Babindamana$^{2}$}
\address{Regis Freguin Babindamana
	\\ Department of Mathematics and Computer Science\\ Universit\'e Marien Ngouabi \\ Facult\'e des Sciences et Techniques \\ BP: 69, Brazzaville, Congo \\
	https://orcid.org/0000-0002-8538-890X}
\email{regis.babindamana@umng.cg}

\author[Bossoto]{Basile Guy Richard Bossoto$^{3}$
	\\
	\\ 
	Universit\'e Marien Ngouabi$^{1,2,3}$\\
	Facult\'e des Sciences et Techniques$^{1,2,3}$\\
	BP: 69, Brazzaville, Congo$^{1,2,3}$\\
	bansimbagilda@gmail.com$^1$, regis.babindamana@umng.cg$^{2}$, basile.bossoto@umng.cg$^3$ \\
}
\address{Basile Guy Richard Bossoto \\ Department of Mathematics\\ Universit\'e Marien Ngouabi \\ Facult\'e des Sciences et Techniques \\ BP: 69, Brazzaville, Congo}
\email{basile.bossoto@umng.cg}

\begin{abstract} \ \\ \rm{
From the results in the literature, the algebraic set of the hyperbola with parameter $n$ defined by $\mathcal{B}_{n}(X, Y, Z)_{\mid_{x\geq 4n}}= \displaystyle \lbrace \left(X: Y: Z\right)\in \mathbb{P}^{2}(\mathbb{Q}) \ \vert \ \displaystyle Y^{2}=X^{2}-4nXZ \rbrace$ where $n$ is a semiprime is proved to be in relation with prime factors of $n$. In the affine space over $\mathbb{Z}_{\geqslant 4n}\times \mathbb{Z}_{\geqslant 0}$, this set has exactly 5 points $\displaystyle\lbrace P_{0}, P_{1}, P_{2}, P_{3}, P_{4}  \rbrace$ with $P_{2}+P_{3}=P_{1}+2P_{2}=P_{4}$ for which knowledge of $P_{2}$ or $P_{3}$ yields the factorization of $n$. However, The non cyclicity of this group structure over rationals and integers and moreover its non good reduction over finite fields constitute the main difficulty in finding its solutions. In this paper we describe an approach to finding $P_{2}$ and $P_{3}$.  We introduce the concept of Hyperbola X-root and Y-root that the solution's greatest common divisors with $n$ reveal prime factors of $n$. We prove that $P_{2}$ and $P_{3}$ can be found on a singular Weierstrass curve isomorphic to a Jacobi quartic using the Hyperbola X-root and Y-root. We present the mathematical framework for this approach.}
\end{abstract}
\maketitle
\section*{Introduction}
In cryptography, the Rivest, Shamir, and Adleman (RSA) Cryptosystem \cite{rsa} security relies on the intractability of factoring large primes.
For cryptanalysis, several factorization methods \cite{fac} have been studied in the literature, mainly the General purpose factorization methods among Fermat's method \cite{fer1, fer2}, the Quadratic Sieve \cite{qs} or the Number Field Sieve \cite{nfs, dv}, and also special purpose methods like the Lentras's method on elliptic curves \cite{ecm, ed}. To date, no computationally efficient factorization algorithm for large semiprimes is known. As an example, recently at the end of 2019 using CADO-NFS \cite{kdo}, the new record of RSA-240, with a semiprime of 795 bits took 953 core years distributed over thousands of computers \cite{loria}. Consequently, this motivates and leads to the exploration of new directions and approaches.
In \cite{xyz}, authors studied the particular hyperbola $\mathcal{B}_{n}(X, Y, Z)_{\mid_{x\geq 4n}}= \displaystyle \lbrace \left(X: Y: Z\right)\in \mathbb{P}^{2}(\mathbb{Q}) \ \vert \ \displaystyle Y^{2}=X^{2}-4nXZ \rbrace$ and presented results among which in affine space, for any RSA semiprime $n$, $\mathcal{B}_{n}(x, y)_{\mid_{\mathbb{Z}}}$ has a generating set: $\mathcal{B}_{n}(x, y)_{\mid_{\mathbb{Z}}}=<P_{0}, P_{1}, P_{2}>$, and $\#\mathcal{B}_{n}(x, y)_{\mid_{\mathbb{Z}}}=3^{2}2^{1}$, and established an equivalence between factoring $n$ and finding one of the $P_{i}'s$, $i=\overline{1,\cdots,3}$. However, this group structure is not cyclic and requires to exploit other approaches to computing the $P_{i}'s$.\\
Towards this, in this paper we present a new approach to finding the $P_{i}^{'}$s thus factoring $n$. We introduce the concept, the Hyperbola X-root and Y-root, denoted respectively Hyper-X and Hyper-Y, such that the solution's the greatest common divisors with $n$ give a prime factor of $n$. We show that $P_{2}$ and $P_{3}$ can be found using a singular Weierstrass curve isomorphic to a Jacobi quartic through the Hyper-X and Hyper-Y, and present some underlying algebraic relations that constitute a new framework introducing parameter relations on some parameters that can be exploited for integer factorization.\\
The paper is organized as follows:\\
- In section one, we present some preliminaries and study some algebraic properties on $\mathcal{B}_{n}(x, y)_{\mid_{\mathbb{Z}}}$ in the affine space, mainly some group classification information, and a discussion on the computation of orders of its subgroups.\\
- In section two, we present the new approach. \\

\textbf{Nomenclature:}\\
Here is a list of the commonly used nomenclature in this paper:\\
$\mathcal{B}_{N}(x, y)_{\mid_{\mathbb{Q}}}=\displaystyle \lbrace \left(x, y\right)\in \mathbb{Q}\times \mathbb{Q} \ / \ y^{2}=x^{2}-4Nx \displaystyle \rbrace$: algebraic set of all rational points on $\mathcal{B}_{N}(x, y)$. \\
$\mathcal{B}_{N}(x, y)_{\mid_{\mathbb{Z}}}=\displaystyle \lbrace \left(x, y\right)\in \mathbb{Z}\times \mathbb{Z} \ / \ y^{2}=x^{2}-4Nx \displaystyle \rbrace$: algebraic set of all integral points on $\mathcal{B}_{N}(x, y)$ \\
$\mathcal{B}_{N}(x, y)_{\mid_{x\geq 4N}}= \lbrace (x, y)\in \mathbb{Z}_{\geq 4N}\times \mathbb{Z}_{\geq 0} \ / \ y^{2}=x^{2}-4Nx\rbrace$: algebraic set of integral points on $\mathcal{B}_{N}(x, y)$ whose $x-$coordinates are greater or equal to $4N$.\\
$Card\displaystyle\left(\mathcal{B}_{N}(x, y)\right)$: the cardinal of $\mathcal{B}_{N}(x, y)$.\\
$\#\Gamma(i)$: the cardinal of $\Gamma(i)$\\
$x | y$: $x$ divides $y$ and $x \nmid y$: $x$ does not divide $y$.\\
Hyper-X: Hyperbola $X-root$ \\
Hyper-Y: Hyperbola $Y-root$

\section{Some Prelimaniries and Algebraic Properties on $\mathcal{B}_{N}(x, y)_{\mid_{\mathbb{Z}}}$} \rm{
In this section, we start by presenting some algebraic properties  on groups and subgroups of $\mathcal{B}_{N}(x, y)_{\mid_{\mathbb{Z}}}$ in the affine space $\mathbb{Z}[x, y]$ underlining the non cyclicity.

\subsection{Some Group Properties}
We know from \cite{xyz} that \\ $Card\displaystyle\left(\mathcal{B}_{N}(x, y)_{\mid_{x\in \mathbb{Z}}} \right)=4\displaystyle\left( Card\displaystyle\left(\mathcal{B}_{N}(x, y)_{\mid_{x\geq 4N}} \right)-1\right) +1_{(4N, 0)}+1_{(0, 0)}=\displaystyle 4\left( Card\displaystyle\left(\mathcal{B}_{N}(x, y)_{\mid_{x\geq 4N}}\right) \right)-2$. \\
And Set $N_{a}$ a prime number, $ \alpha \in \mathbb{Z}_{+}$ such that $N=N_{a}^{\alpha}$. Then $Card\displaystyle\left(\mathcal{B}_{N}(x, y)_{\mid_{x\geq 4N}} \right)=\alpha +1$, and \\  $Card\displaystyle\left(\mathcal{B}_{N}(x, y)_{\mid_{x\in \mathbb{Z}}} \right)=4\alpha +2$ since $Card\displaystyle\left(\mathcal{B}_{N}(x, y)_{\mid_{x\geq 4N}} \right)=\alpha +1$. \\
Considering the case $N=p\times q$, where $p$ and $q$ are primes, then $Card(\mathcal{B}_{N}(x, y)_{\mid_{\mathbb{Z}}})=18$ and we have \\ 
$\mathcal{B}_{N}(x, y)_{\mid_{\mathbb{Z}}}= \displaystyle \lbrace\\
 (p(q + 1)^2, p(q^2 - 1)), \ (p(q +
1)^2, -p(q^2 - 1)), \ (-p(q + 1)^2 + 4n, p(q^2 - 1)), \\ 
(-p(q + 1)^2 +4n, -p(q^2 - 1)), \ (4n, 0), \ (0, 0), \ (q(p + 1)^2, q(p^2 - 1)), \ (q(p +
1)^2, -q(p^2 - 1)), \ (-q(p + 1)^2 + 4n, q(p^2 - 1)), \ (-q(p + 1)^2
+ 4n, -q(p^2 - 1)), \ ((p + q)^2, p^2 - q^2), \ ((p + q)^2, -p^2
+ q^2), \ (-(p + q)^2 + 4n, p^2 - q^2), \ (-(p + q)^2 + 4n, -p^2 + q^2), \ ((n + 1)^2, n^2 - 1), \ ((n + 1)^2, -n^2 +
1), \ (-(n + 1)^2 + 4n, n^2 - 1), \ (-(n + 1)^2 + 4n, -n^2 + 1)  \\
\rbrace$. \\
}
\begin{Proposition} \ \\ \label{prop1} \rm{
$\forall \ N\in \mathbb{Z}_{>0}$, considering the group $\mathcal{B}_{N}(x, y)_{\mid_{\mathbb{Z}}}$ then the following hold:
\begin{itemize}
\item[1)] $\mathcal{B}_{N}(x, y)_{\mid_{\mathbb{Z}}}$ is not a monogenic group, and $\left( \lbrace O=(0, 0), \mathcal{O}=(4N, 0) \rbrace, + \right)$ is an invariant and cyclic subgroup generated by $O=(0, 0)$.
\item[2)] If $N=p\times q$, $p$, $q$ primes, then $\displaystyle\left( \mathcal{B}_{N}(x, y)_{\mid_{\mathbb{Z}}}, \ + \right)$ is generated by the generating part $\lbrace P_{0}, P_{1}, P_{2}\rbrace$.
\end{itemize}
}
\end{Proposition}
\begin{Proof} \ \rm{
\begin{itemize}
\item[1)] Assume $\mathcal{B}_{N}(x, y)_{\mid_{\mathbb{Z}}}$ to be a monogenic group, that's to say, there exists $P\in \mathcal{B}_{N}(x, y)_{\mid_{\mathbb{Z}}}$ such that $<P>= \mathcal{B}_{N}(x, y)_{\mid_{\mathbb{Z}}}$ then $P$ is a $Card\left(\mathcal{B}_{N}(x, y)_{\mid_{\mathbb{Z}}} \right)$-torsion point, ie  $Card\left(\mathcal{B}_{N}(x, y)_{\mid_{\mathbb{Z}}} \right)\cdot P=P_{0}=(4N, 0)$. This yields $P$ to be a $k$-torsion point with $k> 2$ since $Card\left(\mathcal{B}_{N}(x, y)_{\mid_{\mathbb{Z}}}\right) =k=4\alpha+2 \ >2$, which is absurd since $\left( \mathcal{B}_{N}(x, y)_{\mid_{\mathbb{Z}}}\right)_{tors} \cong \mathbb{Z}/2\mathbb{Z}$ which means that $\mathcal{B}_{N}(x, y)_{\mid_{\mathbb{Z}}}$ has only $2$-torsion points. \\
Also $\left( \lbrace O=(0, 0), \mathcal{O}=(4N, 0) \rbrace, + \right)$ is cyclic since $O+O=2\cdot O=2\cdot (0,0)=(4N, 0)$. $\Rightarrow \ <O>=\left( \lbrace O=(0, 0), \mathcal{O}=(4N, 0) \rbrace, + \right)$. The normality of this subgroup is trivial.
\item[2)] $Card\left(\mathcal{B}_{N}(x, y)_{\mid_{x\geq 4N}}\right)=5$ and $\mathcal{B}_{N}(x, y)_{\mid_{x\geq 4N}}=\lbrace P_{0}=(4N, 0), \ P_{1}=((p+q)^2, p^2-q^2), \ P_{2}=(p(q+1)^2, p(q^2-1)), \ P_{3}=(q(p+1)^2, q(p^2-1)), P_{2}=((n+1)^2, n^2-1) \rbrace $. It's already been proved that $P_{3}=P_{2}+P_{1}=P_{2}+P_{1}+P_{0}$ $ \Longrightarrow \ \exists \ i=j=1$ such that $P_{3}=i\cdot P_{2}+j\cdot P_{1}+P_{0}$. And also $P_{4}=P_{2}+P_{3}=2P_{2}+P_{1}+P_{0}$ $\Longrightarrow \ \exists \ i=2, j=1$ such that $P_{4}=i\cdot P_{2}+j\cdot P_{1}+P_{0}$. Therefore $\mathcal{B}_{N}(x, y)_{\mid_{x\geq 4N}}=<P_{0}, P_{1}, P_{2}>$. \\
Now, since all other points over the integers are obtained from symmetry, that's to say, $\forall \ P=(x_{P}, y_{P})\in \mathcal{B}_{N}(x, y)_{\mid_{x\geq 4N}}$, $(x_{P}, y_{P}), (-x_{P}+4N, y_{P}), (-x_{P}+4N, -y_{P})\in \mathcal{B}_{N}(x, y)_{\mid_{\mathbb{Z}}}$, thus $<P_{0}, P_{1}, P_{2}>=\mathcal{B}_{N}(x, y)_{\mid_{\mathbb{Z}}}$.
\end{itemize}
}
\end{Proof}

\subsection{Some properties on Subgroups of $\mathcal{B}_{N}(x, y)_{\mid_{x\in \mathbb{Z}}}$ }
We first recall some important results from the literature, mainly the Cauchy and Lagrange Theorems that we use later.
\begin{Theorem} of Lagrange\ \label{theo} \\ \rm{
If $H$ is a subgroup of a finite group $G$, then the order of $H$ divides the order $G$. That's to say $\vert G \vert = \vert G: H \vert \cdot \vert H \vert$.
}
\end{Theorem}
\begin{Proof} \ \\ \rm{
Can be found in \cite{echo}
}
\end{Proof}

\begin{Theorem} of Cauchy \ \label{cauchy} \\ \rm{ 
Let $G$ be a finite group and $p$ be a prime factor of $|G|$. Then $G$
contains an element of order $p$. Equivalently, $G$ contains a subgroup of order $p$.
}
\end{Theorem}
\begin{Proof} \ \\ \rm{
Can be found in \cite{cauchy1} and \cite{cauchy2} 
}
\end{Proof}
\begin{Lemma} \ \label{lem1} \\ \rm{
$\forall \ N\in \mathbb{Z}_{> 0}$, $Card(\mathcal{B}_{N}(x, y)_{\mid_{\mathbb{Z}}})\equiv 0 \mod 2$, and $\forall \alpha \in \mathbb{Z}_{> 0}$, $p$ and $2\alpha+1$ primes, then there exists a subgroup of $\mathcal{B}_{p^{\alpha}}(x, y)_{\mid_{\mathbb{Z}}} $ of order $2\alpha+1$.
}
\end{Lemma}
\begin{Proof} \ \\ \rm{
This result is trivial and straightforward from \cite{xyz}, since
$Card\displaystyle\left(\mathcal{B}_{N}(x, y)_{\mid_{x\in \mathbb{Z}}} \right)=\displaystyle 2\left( 2Card\displaystyle\left(\mathcal{B}_{N}(x, y)_{\mid_{x\geq 4N}}\right) -1\right)$. \\
Since $\forall \alpha \in \mathbb{Z}_{> 0}$, $Card\left(\mathcal{B}_{p^{\alpha}}(x, y)_{\mid_{x\geq 4p^{\alpha}}}\right)=\alpha+1$, then $Card\displaystyle\left(\mathcal{B}_{p^{\alpha}}(x, y)_{\mid_{x\in \mathbb{Z}}} \right)=4\alpha+2=2(2\alpha+1)$, hence by Theorem \ref{cauchy},  $\mathcal{B}_{p^{\alpha}}(x, y)_{\mid_{\mathbb{Z}}}$ has a subgroup of order $2\alpha+1$.
}
\end{Proof}
\begin{Definition} \ \\ \rm{
	A subgroup of a group $G$ is called a p-Sylow subgroup of $G$ if $p^i$ is the highest power
of the prime p that divides the order of $G$ with $i$ an integer. }
\end{Definition}
\begin{Definition} \ \\ \rm{
For a prime number $p$, the p-adic order or p-adic valuation on $\mathbb{Z}$
is the function 
$v_p: \mathbb{Z}\longrightarrow \mathbb{N}$
\begin{eqnarray*}
	v_{p}(n)=
	\left\lbrace 
	\begin{array}{ll}
		\max \lbrace k \in \mathbb{N}: \ p^k \text{ divides } n \rbrace \ \ \ \ \ \ \text{ if } n \neq 0 \\\
		\infty \ \ \ \ \ \ \text{ if } n = 0
	\end{array} \right.
\end{eqnarray*} }
\end{Definition}
\begin{Example} \ \\ \rm{
	$75=3^1\cdot5^2$ then $2$ is the $5-$adic valuation of $75$ and $1$ is the $3-$adic valuation of $75$. ($v_3 (75)=1$ and $v_5 (75)=2$). }
\end{Example}

\begin{Theorem} \ \\ \rm{
$\forall \ N=p\times q\in \mathbb{Z}_{>0}$, $p$ and $q$ primes, then the following hold:
\begin{itemize}
\item[a)] $\mathcal{B}_{N}(x, y)_{\mid_{x \in \mathbb{Z}}}$ has subgroups of order $2$ and $3$.
\item[b)] $\mathcal{B}_{N}(x, y)_{\mid_{x \in \mathbb{Z}}}$ is a soluble group. 
\item[c)] There exists a 3-Sylow subgroup of $\mathcal{B}_{N}(x, y)_{\mid_{x \in \mathbb{Z}}}$ and $2$ is the 3-adic valuation of the order of $\mathcal{B}_{N}(x, y)_{\mid_{x \in \mathbb{Z}}}$.
\end{itemize}  
}
\end{Theorem}

\begin{Proof} \ \rm{
\begin{itemize}
\item[a)] This proof is a direct result of the Cauchy Theorem
\item[b)] To prove this, we have to prove that there exist $n$, $m$ positive integers, and $p$, $q$ prime numbers such that $Card\displaystyle\left(\mathcal{B}_{N}(x, y)_{\mid_{x\in \mathbb{Z}}} \right)=p^m q^n$. The Proof is straightforward. From \cite{xyz}, if $N=p\times q$, $p$ and $q$ primes, then $Card\displaystyle\left(\mathcal{B}_{N}(x, y)_{\mid_{x\geq 4N}}\right)=5$ and since $Card\displaystyle\left(\mathcal{B}_{N}(x, y)_{\mid_{x\in \mathbb{Z}}} \right)= 4\left( Card\displaystyle\left(\mathcal{B}_{N}(x, y)_{\mid_{x\geq 4N}}\right) \right)-2$, then $Card\displaystyle\left(\mathcal{B}_{N}(x, y)_{\mid_{x\in \mathbb{Z}}} \right)=4\times 5-2=18=2^1 3^2$
\item[c)] From the Theorem of Sylow, since $3$ is prime  and $3 \mid$ the order of $\mathcal{B}_{N}(x, y)_{\mid_{x\in \mathbb{Z}}}$ and from b), $Card\displaystyle\left(\mathcal{B}_{N}(x, y)_{\mid_{x\in \mathbb{Z}}} \right)=2^1 3^2$ and $3\nmid 2$ then $2$ is the 3-adic valuation of $18$.
\end{itemize}}
\end{Proof}

\begin{Remark} \ \\ \rm{
From Lemma \ref{lem1}, $\forall \ \eta\in \mathbb{Z}_{>0}$, $Card(\mathcal{B}_{\eta}(x, y)_{\mid_{\mathbb{Z}}})\equiv 0 \mod 2$. \\
Considering $n=p\times q$ product of two odd primes, then $\nexists \ \eta \in \mathbb{Z}$ such that $Card(\mathcal{B}_{\eta}(x, y)_{\mid_{\mathbb{Z}}})=n$ (From Proposition \ref{prop1}). \\
Now let's consider $\mathcal{B}_{\omega}(x, y)$ such that $Card(\mathcal{B}_{\omega}(x, y)_{\mid_{\mathbb{Z}}})=2n$. Then $Card(\mathcal{B}_{\omega}(x, y)_{\mid_{x\geq 4\omega}})=\frac{n+1}{2}$. \\
Since $n$ is the product of two primes $p$ and $q$, then from Theorem \ref{theo}, $\mathcal{B}_{\omega}(x, y)_{\mid_{\mathbb{Z}}}$ has  subgroups of order $2$, $p$ and $q$. As $\mathcal{B}_{\omega}(x, y)_{\mid_{\mathbb{Z}}}$ is not cyclic, nevertheless by Theorem \ref{cauchy}, it has subgroups of order $2$, $p$ and $q$.
}
\end{Remark}

\subsection{Discussion on the computation of subgroups of order 2, p and q.} \ \rm{
Here we review properties on the order of the group and its subgroups and see how we can compute this.
\begin{itemize}
\item[•] \textbf{Subgroup of order $2$}.\\
From \cite{xyz}, $\displaystyle\left(\lbrace O=(0,0), \ \mathcal{O}=(4\omega, 0)\rbrace, \ + \right)$ is a normal subgroup of $\mathcal{B}_{\omega}(x, y)_{\mid_{\mathbb{Z}}}$,\\
and $<O>=\displaystyle\left(\lbrace O=(0,0), \ \mathcal{O}=(4\omega, 0)\rbrace, \ + \right)$, this shows that $O$ is a point of order $2$. \\
We have $2O=2(0, 0)=\mathcal{O}=(4\omega, 0)$. 
$O$ is a $2-$torsion point, hence an element of order $2$. \\
\item[•] \textbf{Subgroups of order $p$ and $q$}.\\
From the previous subsection, $Card(\mathcal{B}_{\omega}(x, y)_{\mid_{\mathbb{Z}}})=2n$, then $Card(\mathcal{B}_{\omega}(x, y)_{\mid_{x\geq 4\omega}})=\frac{n+1}{2}$. In this case, among the possible values of $\omega$, there exists a prime $p$ such that $\omega=p^{\frac{n-1}{2}}$, particularly for $p=2$. \\
This is almost obvious since from \cite{xyz}, if $p$ is prime and $\alpha \in \mathbb{N}^{\star}$, then $Card(\mathcal{B}_{p^\alpha}(x, y)_{\mid_{x\geq 4p^{\alpha}}})=\alpha+1$. Therefore setting $\alpha+1=\frac{n+1}{2}$, we find $\alpha=\frac{n-1}{2}$. \\
Now it's clear that $Card(\mathcal{B}_{\omega}(x, y)_{\mid_{\mathbb{Z}}})=2n$ with $\omega=2^{\frac{n-1}{2}}$. Therefore, the question of factoring the cardinal $2n$ of 
$\mathcal{B}_{2^{\frac{n-1}{2}}}(x, y)_{\mid_{x\geq 4\cdot 2^{\frac{n-1}{2}}}}$ is an equivalent to finding its subgroups.\\
We are guaranteed from the Theorem of Cauchy (Theorem \ref{cauchy}) that subgroups of $\mathcal{B}_{2^{\frac{n-1}{2}}}(x, y)_{\mid_{x\geq 4\cdot 2^{\frac{n-1}{2}}}}$ have orders, respectively, $2$, $p$, and $q$. \\
\end{itemize}
The problem we have here is that the group is not cyclic, and computing subgroups may not be a practical way to find the orders.\\
For this reason, in the following section we explore a new framework and present the underlying approach }

\section{The new hyperbola based approach to factoring Integers} \ \rm{
In this section, we first present the mathematical framework and then present the approach.\\ 
We consider  $\mathcal{B}_{n}(x, y)_{\mid_{x\geq 4n}}$. 
From arithmetical results on $\mathcal{B}_{n}(x, y)_{\mid_{\mathbb{Z}}}$ group structure in \cite{xyz}, for a RSA cryptosystem modulus, $Card\displaystyle\left(\mathcal{B}_{n}(x, y)_{\mid_{x\geq 4n}}\right)=5$ and $\mathcal{B}_{n}(x, y)_{\mid_{x\geq 4n}}=\lbrace P_{0}=(4n, 0), \ P_{1}=((p+q)^2, p^2-q^2), \ P_{2}=(p(q+1)^2, p(q^2-1)), \ P_{3}=(q(p+1)^2, q(p^2-1)), \ P_{4}=((n+1)^2, n^2-1) \rbrace$ with $\displaystyle\left\langle P_{0}, P_{1}, P_{2} \right\rangle$ as a generating set. It's proved that $P_{4}=P_{3}+P_{2}=P_{1}+2P_{2}$. Another result among others is that factoring $n$ is equivalent to finding any $P \in \mathcal{B}_{n}(x, y)_{\mid_{\mathbb{Z}}} \setminus \displaystyle\left\lbrace O=(0,0), \ P_{0}, \ P_{4}\right\rbrace$. 

\begin{Proposition} \ \\ \rm{
Consider $N=p\times q$, a RSA modulus. Finding integral solutions of $\mathcal{B}_{N}(x, y)_{\mid_{x\geq 4N}}$ thus factoring $N$ is equivalent to solving the system type \begin{eqnarray}
\left\lbrace 
\begin{array}{ll}
X_{1}X_{2}+Y_{1}Y_{2}= 2N(N^{2}+1)\\
X_{2}Y_{1}+X_{1}Y_{2}=2N(N^{2}-1)
\end{array} \right.
\end{eqnarray}
Where $X_{1}, X_{2}, Y_{1}, Y_{2} \ \in \mathbb{Z}_{+}$.}

\end{Proposition}

\begin{Proof} \ \\ \rm{
From \cite{xyz}, $Card\left(\mathcal{B}_{N}(x, y)_{\mid_{x\geq 4N}}\right) = 5$ and $\mathcal{B}_{N}(x, y)_{\mid_{x\geq 4N}}= \lbrace P_{0}=(4N, 0), P_{1}=(x_{P_{1}}, y_{P_{1}}), P_{2}=(x_{P_{2}}, y_{P_{2}}), P_{3}=(x_{P_{3}}, y_{P_{3}}), P_{4}=((N+1)^2, N^2-1)\rbrace$. Since $\lbrace P_{0}, P_{1}, P_{2}  \rbrace$ is a generating part of $\mathcal{B}_{N}(x, y)_{\mid_{x\geq 4N}}$ with respectively $P_{3}=P_{1}+P_{2}$; $P_{4}=P_{2}+P_{3}$, we have from the group law: 
\begin{eqnarray*}
\left\lbrace
\begin{array}{ll}
x_{P_{2}+P_{3}}=\frac{1}{2N}\left[(x_{P_{2}}-2N)(x_{P_{3}}-2N)+y_{P_{2}}y_{  P_{3}}\right]+2N=x_{P_{4}}\\
y_{P_{2}+P_{3}}=\frac{1}{2N}\left[y_{P_{2}}(x_{P_{3}}-2N)+y_{P_{3}}(x_{P_{2}}-2N)\right]=y_{P_{4}}
\end{array} \right. \ \ ;
\left\lbrace 
\begin{array}{ll}
(x_{P_{2}}-2N)(x_{P_{3}}-2N)+y_{P_{2}}y_{P_{3}}= 2N(N^{2}+1)\\
y_{P_{2}}(x_{P_{3}}-2N)+y_{P_{3}}(x_{P_{2}}-2N)=2N(N^{2}-1)
\end{array} \right.
\end{eqnarray*}
Then setting $x_{P_{2}}-2N=X_{1}$, $x_{P_{3}}-2N=X_{2}$, $y_{P_{2}}=Y_{1}$ and $y_{P_{3}}=Y_{2}$, we finally have the system 
\begin{eqnarray*}
\left\lbrace 
\begin{array}{ll}
X_{1}X_{2}+Y_{1}Y_{2}= 2N(N^{2}+1)\\
X_{2}Y_{1}+X_{1}Y_{2}=2N(N^{2}-1)
\end{array} \right.
\end{eqnarray*}
Since the knowledge of $X_{1}, X_{2}, Y_{1}, Y_{2}$ yields to the knowledge of $x_{P_{2}}$, $y_{P_{2}}$, $x_{P_{3}}$, $y_{P_{3}}$ and from \cite{xyz}, $\gcd(x_{P_{2}}, \ N)$ and $\gcd(x_{P_{3}}, \ N)$ yield to non trivial factors of $N$. Hence knowing the solution of $(2)$ is equivalent to factoring $N$.
}
\end{Proof}

\begin{Corollary} \ \\ \rm{
There exists a rational polynomial linking $P_{2}$ to $P_{3}$. 
}
\end{Corollary}

\begin{Proof} \ \\ \rm{
Consider the system $(2)$, after a rearrangement of unknowns, thanks to the commutativity of the multiplication law, the system is
\begin{small}
\begin{eqnarray*}
\left\lbrace 
\begin{array}{ll}
X_{1}X_{2}+Y_{1}Y_{2}= 2N(N^{2}+1)=a\\
Y_{1}X_{2}+X_{1}Y_{2}=2N(N^{2}-1)=b
\end{array} \right. \ \Rightarrow 
\left\lbrace 
\begin{array}{ll}
-X_{1}Y_{2}-Y_{1}\frac{Y_{2}^{2}}{X_{2}}= -a\frac{Y_{2}}{X_{2}}\\
Y_{1}X_{2}+X_{1}Y_{2}=b
\end{array} \right. \ \Rightarrow \ Y_{1}\left(\frac{X_{2}^{2}-Y_{2}^{2}}{X_{2}} \right)=\frac{bX_{2}-aY_{2}}{X_{2}}
\Rightarrow \ \displaystyle Y_{1}=\frac{bX_{2}-aY_{2}}{X_{2}^{2}-Y_{2}^{2}} \\
\text{By the same, } \left\lbrace 
\begin{array}{ll}
-X_{1}\frac{X_{2}^{2}}{Y_{2}}-Y_{1}X_{2}= -a\frac{X_{2}}{Y_{2}}\\
X_{1}Y_{2}+Y_{1}X_{2}=b
\end{array} \right. \ \Rightarrow \ X_{1}\left(\frac{Y_{2}^{2}-X_{2}^{2}}{Y_{2}} \right)=\frac{bY_{2}-aX_{2}}{Y_{2}}
\Leftrightarrow \ \displaystyle X_{1}=\frac{bY_{2}-aX_{2}}{Y_{2}^{2}-X_{2}^{2}} \ \ \ \ \ \ \ \ \ \ \ \ \ \ \ \ \ \ \ \ \ \ \ \
\end{eqnarray*}
\end{small}
 We see that $P_{2}$ can be expressed as a function of $P_{3}$.
}
\end{Proof}
\begin{Definition} \ \\ \rm{
Let $n$ be a RSA cryptosystem modulus and $\displaystyle\mathcal{B}_{n}(x, y)_{\mid_{x\geq 4n}}=\lbrace P_{0}, P_{1}, P_{2}, P_{3}, P_{4} \rbrace$ its related algebraic set .\\
We define the hyperbola X-root and a hyperbola Y-root denoted respectively "Hyper-X" and "Hyper-Y", the values $x_{P_{2}}+x_{P_{3}}$ and $x_{P_{2}}\cdot x_{P_{3}}$ respectively.

}
\end{Definition}

\begin{Lemma} \ \label{lem}\\ \rm{
Consider $n=p\cdot q$, $p$, $q$ primes and $\displaystyle\mathcal{B}_{n}(x, y)_{\mid_{x\geq 4n}}=\lbrace P_{0}=(4n, 0), \ P_{1}=((p+q)^2, p^2-q^2), \ P_{2}=(q(p + 1)^2, q(p^2 - 1)), \ P_{3}=(p(q + 1)^2, p(q^2 - 1)), \ P_{4}=((n+1)^2, n^2-1) \rbrace$.\\
Then the Hyper-X: $X=x_{P_{2}}+x_{P_{3}}$ and Hyper-Y: $Y=x_{P_{2}}\cdot x_{P_{3}}$  verify the equation:
\begin{align*}
R_{n}(X, Y)&=\\
&16n^{4}X^{4}-128n^{5}X^{3}+(320n^{6}-32n^{8}-32n^{4})X^{2}+(128n^{9}-256n^{7}+128n^{5})X-(16n^{6}+16n^{2}-32n^{4})Y^{2}\\
&-64n^{4}X^{2}Y+(384n^{5}+64n^{7}+64n^{3})XY-(192n^{8}+640n^{6}+192n^{4})Y+16n^{12}-64n^{10}+96n^{8}-64n^{6}+16n^{4}\\
&=0
\end{align*}

}
\end{Lemma}
\begin{Proof} \ \\ \label{prof} \rm{
In $\mathcal{B}_{n}(x, y)_{\mid_{x\geq 4n}}$, $P_{1}+2P_{2}=P_{2}+P_{3}=P_{4}$. We then have the system:
\begin{eqnarray*}
\left\lbrace
\begin{array}{ll}
x_{P_{2}+P_{3}}=\frac{1}{2n}\left[(x_{P_{2}}-2n)(x_{P_{3}}-2n)+y_{P_{2}}y_{  P_{3}}\right]+2n=x_{P_{4}} \ \ (a_{1})\\
y_{P_{2}+P_{3}}=\frac{1}{2n}\left[y_{P_{2}}(x_{P_{3}}-2n)+y_{P_{3}}(x_{P_{2}}-2n)\right]=y_{P_{4}} \ \ (a_{2})
\end{array} \right.
\end{eqnarray*}
From $(a_{1})$, $2nx_{P_{4}}-4n^{2}-(x_{p_{1}}-2n)(x_{2p_{2}}-2n)=\displaystyle\left[(x_{1}^{2}-4nx_{1})(x_{2}^{2}-4nx_{2})\right]^{\frac{1}{2}}$. $\Longrightarrow$ \ $(1)\ (x_{1}^{2}-4nx_{1})(x_{2}^{2}-4nx_{2})\geq 0$, \ $(2)\ 2nx_{P_{4}}-4n^{2}\geq(x_{p_{1}}-2n)(x_{2p_{2}}-2n)$ and $(3)\ \displaystyle \left(2nx_{P_{4}}-4n^{2}-(x_{p_{1}}-2n)(x_{2p_{2}}-2n)\right)^{2}=(x_{1}^{2}-4nx_{1})(x_{2}^{2}-4nx_{2})$. $(3)\Longrightarrow -4n^{6}+16n^{5}-8n^{4}x_{p_{1}}-8n^{4}x_{2p_{2}}+4n^{3}x_{p_{1}}x_{2p_{2}}-24n^{4}+16n^{3}x_{p_{1}}-4n^{2}x_{p_{1}}^{2}+16n^{3}x_{2p_{2}}-4n^{2}x_{2p_{2}}^{2}+16n^{3}-8n^{2}x_{p_{1}}-8n^{2}x_{2p_{2}}+4nx_{p_{1}}x_{2p_{2}}-4n^{2}=0$.\\
$\Longleftrightarrow \ -4n^{2}(x_{p_{1}}^{2}+x_{2p_{2}}^{2})-8n^{4}(x_{p_{1}}+x_{2p_{2}})+16n^{3}(x_{p_{1}}+x_{2p_{2}})-8n^{2}(x_{p_{1}}+x_{2p_{2}})+(4n^{3}+4n)x_{p_{1}}x_{2p_{2}}+16n^{5}-4n^{6}-24n^{4}+16n^{3}-4n^{2}=0$ \\
$\Longleftrightarrow \ -4n^{2}(x_{p_{1}}^{2}+x_{2p_{2}}^{2})+(16n^{3}-8n^{4}-8n^{2})(x_{p_{1}}+x_{2p_{2}})+(4n^{3}+4n)x_{p_{1}}x_{2p_{2}}-4n^{6}+16n^{5}-24n^{4}+16n^{3}-4n^{2}=0$
$\Longleftrightarrow \ -4n^{2}(x_{p_{1}}+x_{2p_{2}})^{2}+(16n^{3}-8n^{4}-8n^{2})(x_{p_{1}}+x_{2p_{2}})+(4n^{3}+8n^{2}+4n)x_{p_{1}}x_{2p_{2}}-4n^{6}+16n^{5}-24n^{4}+16n^{3}-4n^{2}=0$
$\Longleftrightarrow \ (x_{p_{1}}+x_{2p_{2}})^{2}+(2n^{2}-4n+2)(x_{p_{1}}+x_{2p_{2}})-(n+\frac{1}{n}+2)x_{p_{1}}x_{2p_{2}}+n^{4}-4n^{3}+6n^{2}-4n+1=0$ \\
$\Longleftrightarrow \ \displaystyle\left(x_{p_{1}}+x_{2p_{2}}+(n-1)^{2}\right)^{2}=\left( n+\frac{1}{n}+2\right)x_{p_{1}}x_{2p_{2}}$. Set $X=x_{p_{1}}+x_{2p_{2}}+(n-1)^{2}$ and $Y=\left( n+\frac{1}{n}+2\right)x_{p_{1}}x_{2p_{2}}$, then we have $Y=X^{2}$ which is a parabola.
\begin{eqnarray}
 \displaystyle \text{We have conditions}
\left\lbrace
\begin{array}{ll}
Y=\left( n+\frac{1}{n}+2\right)x_{p_{1}}x_{2p_{2}} \ \geq 0\\
x_{p_{1}}+x_{p_{2}}=X-(n-1)^{2} \ > 8n
\end{array} \right.
 \ \ \text{since} \ \ 
\left\lbrace
\begin{array}{ll}
x_{p_{1}} \ > 4n\\
x_{p_{2}} \ > 4n
\end{array}  \right.
\Longrightarrow \ x_{p_{1}}+x_{p_{2}} > 8n
\end{eqnarray} 
Set $\mathcal{P}=\lbrace (x, y)\in \mathbb{Z}^{2}/ Y=X^{2}\rbrace$. $\forall (x, y)\in \mathcal{P}$, if $\frac{(n+1)^{2}}{n}\mid Y$ ie $(n+1)^{2}\mid nY$ then 
$x_{p_{1}}+x_{2p_{2}}=X-(n-1)^{2}$ and $x_{p_{1}}x_{2p_{2}}= \frac{Y}{(n+1/n+2)}$.  $n \nmid (n+1)^{2}$, $\Longrightarrow \ n \mid x_{p_{1}}x_{2p_{2}} \ \Longrightarrow \ \exists k_{1}\in \mathbb{Z}\ / \ x_{p_{1}}x_{2p_{2}}=k_{1}n \ \Longrightarrow \ Y=\frac{(n+1)^{2}}{n}k_{1}n=k_{1}(n+1)^{2} \ / \ k_{1}=k^{2}$, $k\in \mathbb{Z}$. Since $Y=(x_{p_{1}}+x_{2p_{2}})^{2}=(n+1)^{2}k^{2}$, $\Longrightarrow \ X=\pm k(n+1)$ and taking into account the conditions in $(2)$ we consider the case  $X=k(n+1)$. \\
$Y=X^{2} \ \Longleftrightarrow \ \displaystyle\left[(x_{p_{1}}+x_{2p_{2}})+(n-1)^{2}\right]^{2}=\frac{(n+1)^{2}}{n}k^{2}n \ \Longrightarrow \ x_{p_{1}}+x_{2p_{2}}=k(n+1)-(n-1)^{2}$ and $x_{p_{1}}x_{2p_{2}}=k^{2}n$. \\


hence from $(a_{1})$, $\exists \ k\in \mathbb{Z}_{>0}$ such that $x_{P_{2}}+x_{P_{3}}=k(n+1)-(n-1)^2$ and $x_{P_{2}}x_{P_{3}}=k^2n$, with the relation $(x_{P_{2}}+x_{P_{3}}+(n-1)^2)^2=(n+\frac{1}{n}+2)x_{P_{2}}x_{P_{3}}$. \\
Let's consider now $(a_{2})$:
\begin{small}
\begin{align*}
y_{P_{2}}(x_{P_{3}}-2n)+y_{P_{3}}(x_{P_{2}}-2n)=2ny_{n} \ & \Longleftrightarrow \ (x_{P_{3}}-2n)\sqrt{x_{P_{2}}^{2}-4nx_{P_{2}}}+(x_{P_{2}}-2n)\sqrt{x_{P_{3}}^{2}-4nx_{P_{3}}}=2ny_{n}\\
\left[(x_{P_{3}}-2n)\sqrt{x_{P_{2}}^{2}-4nx_{P_{2}}}-2ny_{n} \right]&=-(x_{P_{2}}-2n)\sqrt{x_{P_{3}}^{2}-4nx_{P_{3}}} \\
 \Longleftrightarrow \ \left[2ny_{n}-(x_{P_{3}}-2n)\sqrt{x_{P_{2}}^{2}-4nx_{P_{2}}} \right]&=(x_{P_{2}}-2n)\sqrt{x_{P_{3}}^{2}-4nx_{P_{3}}} \\
  \Longleftrightarrow \ \left[2ny_{n}-(x_{P_{3}}-2n)\sqrt{x_{P_{2}}^{2}-4nx_{P_{2}}} \right]^{2}&=(x_{P_{2}}-2n)^{2}(x_{P_{3}}^{2}-4nx_{P_{3}}) \\
\Longleftrightarrow \ 4n^{2}y_{n}^{2}+(x_{P_{3}}-2n)^{2}(x_{P_{2}}^{2}-4nx_{P_{2}})-4&ny_{n}(x_{P_{3}}-2n)\sqrt{x_{P_{2}}^{2}-4nx_{P_{2}}} = (x_{P_{2}}-2n)^{2}(x_{P_{3}}^{2}-4nx_{P_{3}}) \\ 
\Longleftrightarrow \ [4n^{2}y_{n}^{2}+(x_{P_{3}}-2n)^{2}(x_{P_{2}}^{2}-4nx_{P_{2}})  - & (x_{P_{2}}-2n)^{2}(x_{P_{3}}^{2}-4nx_{P_{3}})]=4ny_{n}(x_{P_{3}}-2n)\sqrt{x_{P_{2}}^{2}-4nx_{P_{2}}} \\
\Longleftrightarrow \ [4n^{2}y_{n}^{2}+(x_{P_{3}}-2n)^{2}(x_{P_{2}}^{2}-4nx_{P_{2}}) - & (x_{P_{2}}-2n)^{2}(x_{P_{3}}^{2}-4nx_{P_{3}})]^{2}=16n^{2}y_{n}^{2}(x_{P_{3}}-2n)^{2}(x_{P_{2}}^{2}-4nx_{P_{2}})
\end{align*}
\begin{align*} 
&16n^{12} - 64n^{10} + 96n^8 - 128n^5(xp_{2} + xp_{3})^3 + 16n^4(xp_{2} +
xp_{3})^4 - 64n^4(xp_{2} + xp_{3})^2\cdot xp_{2}\cdot xp_{3} - 64n^6\\ 
&- 16(n^6 - 2n^4 +n^2)xp_{2}^2\cdot xp_{3}^2 + 16n^4 + 64(n^7 + 6n^5 + n^3)(xp_{2} + xp_{3})xp_{2}\cdot xp_{3} -32(n^8 - 10n^6 + n^4)(xp_{2}
+xp_{3})^2 \\
&-64(3n^8 + 10n^6 +3n^4)xp_{2}\cdot xp_{3} + 128(n^9 - 2n^7 + n^5)(xp_{2} + xp_{3})16 \, n^{12} - 64 \, n^{10} + 96 \, n^{8} - 128 \,
{\left(\mathit{xp}_{2}^{3} + \mathit{xp}_{3}^{3}\right)} n^{5}\\
& - 64 \,n^{6} + 16 \, {\left(\mathit{xp}_{2}^{4} + \mathit{xp}_{3}^{4}\right)}
n^{4} - 16 \, {\left(n^{6} + n^{2}\right)} \mathit{xp}_{2}^{2}
\mathit{xp}_{3}^{2}+ 16 \, n^{4} + 64 \, {\left(n^{7} + n^{3}\right)}
{\left(\mathit{xp}_{2} + \mathit{xp}_{3}\right)}\mathit{xp}_{2}\mathit{xp}_{3} \\
&- 256 \, {\left(n^{8} + n^{4}\right)}\mathit{xp}_{2}\mathit{xp}_{3} - 32 {\left(n^{8} - 10 \, n^{6} + n^{4}\right)}{\left(\mathit{xp}_{2}^{2}  + \mathit{xp}_{3}^{2}\right)}
 + 128 \,{\left(n^{9} - 2 \, n^{7} + n^{5}\right)} {\left(\mathit{xp}_{2} +
\mathit{xp}_{3}\right)}=0
\end{align*}
\end{small}
Setting the hyper-X $xp_{2}+xp_{3}=X$ and Hyper-Y  $xp_{2}\cdot xp_{3}=Y$, we finally have \\  
$R_{n}(X, Y)= 16n^{4}X^{4}-128n^{5}X^{3}+(320n^{6}-32n^{8}-32n^{4})X^{2}+(128n^{9}-256n^{7}+128n^{5})X-(16n^{6}+16n^{2}-32n^{4})Y^{2}-64n^{4}X^{2}Y+(384n^{5}+64n^{7}+64n^{3})XY-(192n^{8}+640n^{6}+192n^{4})Y+16n^{12}-64n^{10}+96n^{8}-64n^{6}+16n^{4}=0$.
}
\end{Proof}
\begin{Proposition} \ \\ \label{prr}\rm{
$k$ as defined in the above Proof (Proof \ref{prof}), is of the form $k=2(n-\varepsilon)$, $\varepsilon \in \mathbb{Z}_{>0}$. \\
And the Hyper-X $X=2(n-\varepsilon)(n+1)-(n-1)^2$ and Hyper-Y $Y=4n(n-\varepsilon)^2$.
}
\end{Proposition}
\begin{Proof} \ \\ \rm{
From $(1)$ , $k=(p+1)(q+1)=n+(p+q)+1$. Since $p$ and $q$ are odd primes, then $p+q$ is even and $p+q+1$ is odd. Now, since $n$ is an odd semiprime and $p+q+1$  odd, then $k$ is even. In this case, $ \exists \ k^{'}\in \mathbb{Z}_{>0}$ such that $k=2k^{'}$ with $k^{'}=\frac{n+(p+q)+1}{2}$.\\
Let's now evaluate the sign of $n-k^{'}$.\\
$n-k^{'}=n-\frac{n+(p+q)+1}{2}=\frac{n-(p+q)-1}{2}>0$ since $pq>p+q+1$. Then $n>k^{'}$. Hence $ \exists \ \varepsilon \in \mathbb{Z}_{>0}$ such that $k^{'}=n-\varepsilon$, finally $k=2(n-\varepsilon)$. \\
Furthermore from $(1)$, the Hyper-X and Hyper-Y, respectively, can be expressed in terms of $k$ by $X=x_{P_{2}}+x_{P_{3}}=k(n+1)-(n-1)^2$ respectively $Y=x_{P_{2}}x_{P_{3}}=k^2n$. \\
Now plugging $k=2(n-\varepsilon)$ into the Hyper-X and Hyper-Y, we get the Hyper-X: $X=2(n-\varepsilon)(n+1)-(n-1)^2$ and Hyper-Y: $Y=4n(n-\varepsilon)^2$.
}
\end{Proof}

\subsection{Approaches Computing of $\varepsilon$}
In this subsection, we present results on approaches to computing $\varepsilon$.\\
First, let's consider the following relations:
\begin{Proposition} (parameters in terms of $\varepsilon$) \ \\ \rm{
Consider $n_{\varepsilon}=p_{\varepsilon}q_{\varepsilon}$ and $\varphi_{\varepsilon}$ the Euler totient function of $n$. Then we have: 
\begin{enumerate}
\item $p_{\varepsilon}+q_{\varepsilon}=n-1-2\varepsilon$
\item $p_{\varepsilon}-q_{\varepsilon}=\sqrt{(n-1-2\varepsilon)^2-4n}$
\item $n_{\varepsilon}=2\varepsilon+1+p_{\varepsilon}+q_{\varepsilon}$
\item $\varphi_{\varepsilon}=2(\varepsilon+1)$
\item $X=2(n-\varepsilon)(n+1)-(n-1)^2$
\item $X=n^2+4n-2\varepsilon(n+1)-1$
\item $X=n^2-(2\varepsilon-4)(n+1)-5$
\item More generally the hyper-X is of the form $n^2-an-b$, $ a, b > 0$ with $b-a=5$.
\end{enumerate}
}
\end{Proposition}
\begin{Proof} \ \rm{
\begin{enumerate}
\item Considering the Hyper-Y $Y=4n(n-\varepsilon)^2$. Since $Y=x_{P_{2}}x_{P_{3}}=n(p+1)^2(q+1)^2=n(n+p+q+1)^2$, then $n+(p_{\varepsilon}+q_{\varepsilon})+1=2(n-\varepsilon)$, hence $p_{\varepsilon}+q_{\varepsilon}=n-1-2\varepsilon$.
\item $p_{\varepsilon}+q_{\varepsilon}$ and $n$ satisfy $T^2-(p_{\varepsilon}+q_{\varepsilon})T+n_{\varepsilon}=0$. In this case after a simple computation of roots, it outcomes that $p_{\varepsilon}=\frac{n-1-2\varepsilon+\sqrt{(n-1-2\varepsilon)^2-4n}}{2}$ and $q_{\varepsilon}=\frac{n-1-2\varepsilon-\sqrt{(n-1-2\varepsilon)^2-4n}}{2}$. Hence $p_{\varepsilon}-q_{\varepsilon}=\sqrt{(n-1-2\varepsilon)^2-4n}$
\item $p_{\varepsilon}+q_{\varepsilon}=n_{\varepsilon}-1-2\varepsilon$, then $n_{\varepsilon}=2\varepsilon+1+p_{\varepsilon}+q_{\varepsilon}$.
\item $\varphi=n-(p+q)+1$ then $\varphi_{\varepsilon}=n-(p_{\varepsilon}+q_{\varepsilon})+1=n+1-n+1+2\varepsilon=2(\varepsilon+1)$.
\item $X=2(n-\varepsilon)(n+1)-(n-1)^2$
is straightforward since from $(1)$, $k=2(n-\varepsilon)$.
\item By expanding $(6)$, $X=n^2+4n-2\varepsilon(n+1)-1$ is trivial. 
\item $X=n^2+4n-2\varepsilon(n+1)-1=n^2+4n-2\varepsilon(n+1)-1-4+4=n^2+4(n+1)-2\varepsilon(n+1)-1-4$ Then $X=n^2-(2\varepsilon-4)(n+1)-5$.
\item This results from an observation from $(7)$. Indeed $X=n^2-(2\varepsilon-4)(n+1)-5=n^2-(2\varepsilon-4)n-2\varepsilon+4-5=n^2-(2\varepsilon-4)n-(2\varepsilon+1)$, then $a=2\varepsilon-4$ and $b=2\varepsilon+1$, with $b-a=5$. It's now obvious that $a+k=2n-4$ considering the Proposition \ref{prr}.
\end{enumerate}
}
\end{Proof}

\begin{Conjecture} \ \ \\ \rm{
Let $f_{1}(\varepsilon)=4\varepsilon^2+4(1-n)\varepsilon+n^2-6n+1$, 
$f_{2}(\varepsilon)=4(n^2-2n+1)\varepsilon^2+4(3n^2-n^3-3n+1)\varepsilon+n^4-8n^3+14n^2-8n+1$ and $\Gamma(i)=\lbrace i\in \mathbb{Z}_{>0}: f_{1}(i)=y^2 \text{ or } f_{2}(i)=z^2 \rbrace$.\\ Then
$\#\Gamma(i)=3$ and $\Gamma(i)=\lbrace\varepsilon_{1}, \varepsilon_{2}, n \rbrace$ with $\varepsilon_{1}+\varepsilon_{2}=n-1$ and $f_{2}(\varepsilon_{1})=f_{1}(n)f_{1}(\varepsilon_{1})$. In this case, $\varepsilon=\displaystyle\min_{i\in\mathbb{Z}_{>0}} \Gamma(i)=\min\lbrace \varepsilon_{1}, \varepsilon_{2} \rbrace$ and $\displaystyle\max_{i\in\mathbb{Z}_{>0}} \Gamma(i)=n$.
}
\end{Conjecture}
\begin{Lemma} \ \rm{
\begin{enumerate}
\item $\varphi_{\varepsilon_{1}}(n)+\varphi_{\varepsilon_{2}}(n)=2(n+1)$ and $\varphi_{\varepsilon_{1}}(n)\varphi_{\varepsilon_{2}}(n)=4(\varepsilon_{1}\varepsilon_{2}+n)$.
\item $X_{\varepsilon_{1}}+X_{\varepsilon_{2}}=8n$ and $X_{\varepsilon_{1}}-X_{\varepsilon_{2}}=4k^{'}$
\end{enumerate}
}
\end{Lemma}
\begin{Proof} \ \rm{
\begin{enumerate}
\item $\varphi_{\varepsilon_{1}}(n)+\varphi_{\varepsilon_{2}}(n)=2(n+1)=2(\varepsilon_{1}+1)+2(\varepsilon_{1}+1)=2(\varepsilon_{1}+\varepsilon_{2}+2)$.\\ From the previous Conjecture, $\varepsilon_{1}+\varepsilon_{2}=n-1$, then $\varphi_{\varepsilon_{1}}(n)+\varphi_{\varepsilon_{2}}(n)=2(n-1+2)=2(n+1)$. \\
Also $\varphi_{\varepsilon_{1}}(n)\varphi_{\varepsilon_{2}}(n)=4(\varepsilon_{1}+1)(\varepsilon_{2}+1)=4(\varepsilon_{1}\varepsilon_{2}+\varepsilon_{1}+\varepsilon_{2}+1)=4(\varepsilon_{1}\varepsilon_{2}+n-1+1)=4(\varepsilon_{1}\varepsilon_{2}+n)$.

\item $X_{\varepsilon_{1}}+X_{\varepsilon_{2}}=n^2+4n-2\varepsilon_{1}(n+1)-1+n^2+4n-2\varepsilon_{2}(n+1)-1=2n^{2}+8n-2n(\varepsilon_{1}+\varepsilon_{2})-2(\varepsilon_{1}+\varepsilon_{2})-2$\\
$2n^2+8n-(n-1)-2(n-1)-2=2n^2+8n-2n^2+2n-2n+2-2=8n$\\
Also $X_{\varepsilon_{1}}-X_{\varepsilon_{2}}=n^2+4n-2\varepsilon_{1}(n+1)-1-n^2-4n+2\varepsilon_{1}(n+1)+1=2(n+1)(\varepsilon_{2}-\varepsilon_{1})=4k^{'}$ with $k^{'}=\frac{n+1}{2}(\varepsilon_{2}-\varepsilon_{1})$, since $n$ is an odd semiprime and $n+1$ is even.

\end{enumerate}

}
\end{Proof}

\begin{Remark} \ \ \label{remark} \rm{
\begin{itemize}
	\item[$(1)$] consider the map $h: \lbrace P_{2}, P_{3}\rbrace \subset \mathcal{B}_{n}(x, y)_{\mid_{x\geq 4n}} \longrightarrow \ R_{n}(X, Y)$, $\lbrace P_{2}, P_{3}\rbrace\longmapsto \  \left(x_{P_{2}}+x_{P_{3}}, \ x_{P_{2}}x_{P_{3}}\right)$, it's trivial that $\ker (h)=\lbrace P_{2}, P_{3}\rbrace$. 
	\item[$(2)$] $R_{n}(X, Y)=R_{X}+R_{XY}+\displaystyle\gamma$ verifying:\\
	$R_{X}(4n)=R_{X}(0)=0$\\
	$R_{X}((n+1)^2)+\displaystyle\gamma = R_{X}(-(n-1)^2)+\displaystyle\gamma = 0$\\
	$R_{XY}(4n, 4n^2)=0$.\\ 
	with $R_{X}=16n^{4}X^{4}-128n^{5}X^{3}+(320n^{6}-32n^{8}-32n^{4})X^{2}+(128n^{9}-256n^{7}+128n^{5})X$; \\ $R_{XY}=-(16n^{6}+16n^{2}-32n^{4})Y^{2}-64n^{4}X^{2}Y+(384n^{5}+64n^{7}+64n^{3})XY-(192n^{8}+640n^{6}+192n^{4})Y$; \\ and $\displaystyle\gamma =16n^{12}-64n^{10}+96n^{8}-64n^{6}+16n^{4}$
\end{itemize}
}
\end{Remark}
\begin{Theorem} \ \\ \ \label{theoB} \rm{
 The equation $R_{n}(X, Y)$ is a Jacobi quartic given by $J_{Q}(X, Z): \ Z^{2}=aX^{4}+bX^{3}+cX^{2}+dX+e^2$\\ where $a=1$; \ $b=-16n$;  $c=2n^4+92n^2+2$; \ $d=-(16n^5+224n^3+16n)$ and \ $e^2=n^8+28n^6+198n^4+28n^2+1$ \\  through the isomorphism $\mathcal{I}: \  R_{n}(X, Y) \longrightarrow \ J_{Q}(X, Z)$,\\ 
$\left(X, Y \right) \longmapsto$ $\left(X, \ \ 	\frac{(n^2-1)^2}{n(n^2+1)}\left( Y+\frac{64n^{4}X^{2}-(384n^{5}+64n^{7}+64n^{3})X+(192n^{8}+640n^{6}+192n^{4})}{32(n^3-n)^2}\right) \right)$.

\begin{Proof} \ \ \\ \rm{
	Setting $R_{n}(X, Y)=R_{X}+R_{XY}+\displaystyle\gamma$, 
\begin{align*}
R_{XY} &= -(16n^{6}+16n^{2}-32n^{4})Y^{2}-64n^{4}X^{2}Y+(384n^{5}+64n^{7}+64n^{3})XY-(192n^{8}+640n^{6}+192n^{4})Y\\
	&= -(16n^{6}+16n^{2}-32n^{4})Y^{2}-\left[64n^{4}X^{2}-(384n^{5}+64n^{7}+64n^{3})X+(192n^{8}+640n^{6}+192n^{4}) \right]Y\\
	&= -(16n^{6}+16n^{2}-32n^{4})\displaystyle\left( Y^{2}+\frac{\left[64n^{4}X^{2}-(384n^{5}+64n^{7}+64n^{3})X+(192n^{8}+640n^{6}+192n^{4}) \right]}{(16n^{6}+16n^{2}-32n^{4})}Y\right) \\
	&= -(16n^{6}+16n^{2}-32n^{4})\displaystyle\left(\left[  Y+\frac{64n^{4}X^{2}-(384n^{5}+64n^{7}+64n^{3})X+(192n^{8}+640n^{6}+192n^{4})}{2(16n^{6}+16n^{2}-32n^{4})}\right] ^{2}\right)\\ 
	& \ \ \ \ \ \ \ -(16n^{6}+16n^{2}-32n^{4})\displaystyle\left( -\frac{\left[64n^{4}X^{2}-(384n^{5}+64n^{7}+64n^{3})X+(192n^{8}+640n^{6}+192n^{4})\right]^{2}}{4(16n^{6}+16n^{2}-32n^{4})^2}\right) \\
	&=  -\left[4(n^3-n)\right]^2\displaystyle\left(\left[  Y+\frac{64n^{4}X^{2}-(384n^{5}+64n^{7}+64n^{3})X+(192n^{8}+640n^{6}+192n^{4})}{2\left[4(n^3-n)\right]^2}\right] ^{2}\right)\\
	& \ \ \ \ \ \ \ \ \ \ \ \ \ \ \ \ \ \ \ \ \ \ \ \ \ \ \  +\frac{\left[64n^{4}X^{2}-(384n^{5}+64n^{7}+64n^{3})X+(192n^{8}+640n^{6}+192n^{4})\right]^{2}}{4\left[4(n^3-n) \right]^2} \\
	&= -Z_{\star}^2+B_{X}^2 \\
	& \text{  With } \ Z_{\star}^2=\left[4(n^3-n)\right]^2\displaystyle\left[  Y+\frac{64n^{4}X^{2}-(384n^{5}+64n^{7}+64n^{3})X+(192n^{8}+640n^{6}+192n^{4})}{2\left[4(n^3-n)\right]^2}\right] ^{2} \\
	& \text{  and } B_{X}^2=\frac{\left[64n^{4}X^{2}-(384n^{5}+64n^{7}+64n^{3})X+(192n^{8}+640n^{6}+192n^{4})\right]^{2}}{4\left[4(n^3-n) \right]^2} \\
	& \text{ Now } R_{n}(X,Y)=R_{X}+R_{XY}+\displaystyle\gamma=R_{X}-Z_{\star}^2+B_{X}^2+\displaystyle\gamma=0. \text{ \ \ Then } Z_{\star}^2=R_{X}+B_{X}^2+\displaystyle\gamma \\
	Z_{\star}^2&=R_{X}+B_{X}^2+\displaystyle\gamma \\
	&=16n^{4}X^{4}-128n^{5}X^{3}+(320n^{6}-32n^{8}-32n^{4})X^{2}+(128n^{9}-256n^{7}+128n^{5})X+\\
	&\frac{\left[64n^{4}X^{2}-(384n^{5}+64n^{7}+64n^{3})X+(192n^{8}+640n^{6}+192n^{4})\right]^{2}}{4\left[4(n^3-n) \right]^2} +16n^{12}-64n^{10}+96n^{8}-64n^{6}+16n^{4}\\
	&=\frac{16(X^2-8nX+n^4+14n^2+1)(n^2+1)^2n^4}{(n+1)^2(n-1)^2}\\
	&=\frac{16n^4(n^2+1)^2}{(n^2-1)^2}\left[X^2-8nX+n^4+14n^2+1\right]^2\\
	&=\frac{16n^4(n^2+1)^2}{(n^2-1)^2}\left[X^4-16nX^3+(2n^4+92n^2+2)X^2-(16n^5+224n^3+16n)X+n^8+28n^6+198n^4+28n^2+1\right]\\
	& \text{Setting } Z^2=\frac{(n^2-1)^2}{16n^4(n^2+1)^2} Z_{\star}^2 \text{, then we have :}\\
	Z^2&=X^4-16nX^3+(2n^4+92n^2+2)X^2-(16n^5+224n^3+16n)X+n^8+28n^6+198n^4+28n^2+1 \\
	& \text{By identification: }\\
	& \text{$a=1$; \ $b=-16n$; \ $c=2n^4+92n^2+2$; \ $d=-(16n^5+224n^3+16n)$; \ $e^2=n^8+28n^6+198n^4+28n^2+1$}
\end{align*}
}
\end{Proof}
}
\end{Theorem}


\begin{Remark} \ \\ \rm{
$X$ and $Y$ defined as in Lemma \ref{lem}, verify the parabolic equation:
\begin{small}
\begin{align*}
	Y=\displaystyle\left[\frac{n(n^2+1)}{(n^2-1)^2}-\frac{2n^4}{n^6-2n^4+n^2} \right]X^2+ \displaystyle\left[\frac{2(n^7+6n^5+n^3)}{n^6-2n^4+n^2}-\frac{8n^2}{(n^2-1)^2} \right]X+\displaystyle\frac{(n^4+14n^2+1)(n^2+1)n}{(n^2-1)^2}-\frac{6n^8+20n^6+6n^4}{n^6-2n^4+n^2} 
\end{align*}
\end{small}

}
\end{Remark}

\textbf{Discussion:}\\
 Considering  $J_{Q}(X, Z): \ Z^{2}=aX^{4}+bX^{3}+cX^{2}+dX+e^2$, we are interested in its rational solutions. In the proof of Theorem \ref{theoB}, we rearranged the expression by affecting the common rational term to $Z_{\star}^2$ to get rid of the denominators in the right expression, thus having only integral coefficients $a$, $b$, $c$, $d$ and $e^2$ on the right side of the equation. \\ This trick is important to apply the Rational Root Theorem \cite{rat1} for $Z^{2}=0$. Now being under the applicability conditions of Rational Root Theorem, we can state that each rational solution $X = p_{i}⁄q_{i}$, written in the lowest terms so that $p_{i}$ and $q_{i}$ are relatively prime, satisfies both \textit{$(1)$} $p_{i}$ is an integer factor of the constant term $e^{2}$, and \textit{$(2)$} $q_{i}$ is an integer factor of the leading coefficient $a=1$.\\
We get for $X=n^{4} + 14n^{2} + 1$, \ $Z=\pm(n^4 + 14n^2 - 8n + 2)(n^4 + 14n^2 + 1)$.\\
 The remaining question is how we can compute the points of this quartic to efficiently find its integral solutions. 
 To this question, we try to find an isomorphic structure with good arithmetic properties like having a group law.\\
 Now, before exposing the isomorphic structure, let us first consider the following result.
 
 \begin{Proposition} \ \\ \rm{
Given $(X, Z)\in J_{Q}(X, Z)$, $(X, Z)$ yields to the sought solution of $(X, Y)\in R_{n}(X, Y)$ if and only if $(X, Y)\in \mathbb{Z}_{+}^{2}$.
} 	
 \end{Proposition}

\begin{Proof} \ \\ \rm{
Here we give the shortest proof that is straightforward, since $Z^{2}=aX^{4}+bX^{3}+cX^{2}+dX+e^2$ with integer coefficients, and from Lemma \ref{lem}, $X$ solution then $X\in \mathbb{Z}$, then $Z^2$ is an integer for the $X$ solution. Hence $Z\in \mathbb{Z}$. 
}	
\end{Proof}

\begin{Proposition} \ \\ \rm{
The Jacobi quartic $J_{Q}(X, Z)$ is singular.	
}
\end{Proposition}
\begin{Proof} \ \\ \rm{
Let's consider the right polynomial in $X$. Set $k(X)=aX^{4}+bX^{3}+cX^{2}+dX+e^2$ and $k^{'}(X)=4aX^3+3bX^2+2cX+d$.\\
We know that if $D$ is a discriminant of $k(x)$, then we have the relation $R(k, k^{'})=(-1)^{\frac{n(n-1)}{2}}a_{n}D$ where $R$ is the resultant. \\
Let's use the Sylvester's $(4+3)\times (4+3)$ matrix to calculate the resultant of $k$ and $k^{'}$.\\ \\
$R(k, k^{'})=
\begin{vmatrix}
	a & b & c & d & e & 0 & 0 \\
	0 & a & b & c & d & e & 0 \\
	0 & 0 & a & b & c & d & e \\
	4a & 3b & 2c & d & 0 & 0 & 0 \\
	0 & 4a & 3b & 2c & d & 0 & 0 \\
	0 & 0 & 4a & 3b & 2c & d & 0 \\
	0 & 0 & 0 & 4a & 3b & 2c & d
\end{vmatrix}$ \\ \\
$=a b^{2} c^{2} d^{2} - 4 \, a^{2} c^{3} d^{2} - 4 \, a b^{3} d^{3} + 18 \, a^{2} b c d^{3} - 27 \, a^{3} d^{4} - 4 \, a b^{2} c^{3} e + 16 \, a^{2} c^{4} e + 18 \, a b^{3} c d e - 80 \, a^{2} b c^{2} d e - 6 \, a^{2} b^{2} d^{2} e + 144 \, a^{3} c d^{2} e - 27 \, a b^{4} e^{2} + 144 \, a^{2} b^{2} c e^{2} - 128 \, a^{3} c^{2} e^{2} - 192 \, a^{3} b d e^{2} + 256 \, a^{4} e^{3}
$\\
Then $D=b^{2} c^{2} d^{2} - 4 \, a c^{3} d^{2} - 4 \, b^{3} d^{3} + 18 \, a b c d^{3} - 27 \, a^{2} d^{4} - 4 \, b^{2} c^{3} e + 16 \, a c^{4} e + 18 \, b^{3} c d e - 80 \, a b c^{2} d e - 6 \, a b^{2} d^{2} e + 144 \, a^{2} c d^{2} e - 27 \, b^{4} e^{2} + 144 \, a b^{2} c e^{2} - 128 \, a^{2} c^{2} e^{2} - 192 \, a^{2} b d e^{2} + 256 \, a^{3} e^{3}$, and considering $a=1$; \ $b=-16n$; \ $c=2n^4+92n^2+2$; \ $d=-(16n^5+224n^3+16n)$; \ $e^2=n^8+28n^6+198n^4+28n^2+1$, $D=0$. Therefore $k(x)$ admits points of multiple order. Hence $J_{Q}(X, Z)$ is singular.
	

}
\end{Proof}

\begin{Remark} \ \\ \rm{
Considering $J_{Q}(X, Z)$, in our context of Hyper-X and Hyper-Y, finding special points are important  since they should lead to locate rational points on a new structure, here Weierstrass singular curves through isomorphism. We present some special points:\\
\begin{small}
$\displaystyle\left(4n, \ \frac{2n^5-2n^3+(n^2+1)(n^2-1)n^2}{(n^2-1)n} \right) =(X, Y)$ then $\displaystyle\left(4n, ( n^2-1)^2\right) =(X, Z)$; \\
$\left( (n+1)^2, \ 0\right)=(X, Y)$ then $\displaystyle\left((n+1)^2, \ -2(n^2+1)(n-1)^2 \right)=(X, Z)$;  \\
$\displaystyle\left(n^{4} + 14n^{2} + 1, \frac{n(n^4 + 15n^2 - 2n + 2)^{2}}{(n+1)^2}\right)=(X, Y)$ then $\displaystyle\left(n^{4} + 14n^{2} + 1, \ \ (n^4 + 14n^2 - 8n + 2)(n^4 + 14n^2 + 1)\right)=(X, Z)$\\
$\displaystyle\left(-(n^{4} + 14n^{2} + 1), \frac{n(n^4 + 15n^2 - 2n + 2)^{2}}{(n+1)^2}\right)=(X, Y)$ then $\displaystyle\left(-(n^{4} + 14n^{2} + 1), \ \ (n^4 + 14n^2 + 8n + 2)(n^4 + 14n^2 + 1)\right)=(X, Z)$ 
\end{small}
\begin{Example} \ \ \\
if $n=15$, $\displaystyle\left(n^2-2n-7, \ 4n(n - 3)^2\right)=(X, Y)$ then $\displaystyle\left(n^2-2n-7, \ \ 2(n^2 - 8n + 25)(n + 1)^2\right)=(X, Z)$\\
if $n=21$, $\displaystyle\left(n^2-6n-11, \ 4n(n - 5)^2 \right)=(X, Y)$ then $\displaystyle\left(n^2-6n-11, \ \ 2(n^2 - 12n + 61)(n + 1)^2\right)=(X, Z)$\\
if $n=35$, $\displaystyle\left(n^2-18n-23, \ 4n(n - 11)^2 \right)=(X, Y)$ then $\displaystyle\left(n^2-18n-23,+ \ \ 2(n^2 - 24n + 265)(n + 1)^2\right)=(X, Z)$
\end{Example}

}	
\end{Remark}


\begin{Theorem} \ \\ \ \rm{
The equation $J_{Q}(X, Z)$ is isomorphic to the Weierstrass equation $E_{w}(x, y)_{\mid_{\mathbb{Q}}}: \ y^2+a_{1}xy+a_{3}y=x^3+a_{2}x^2+a_{4}x+a_{6}$\\
with $a_{1}=\displaystyle\frac{d}{e}$, \ $a_{3}=2eb$, \ $a_{2}=\displaystyle c-\frac{d^2}{4e^2}$, \ $a_{4}=-4e^2a$, \ $a_{6}=\displaystyle a\left(d^2-4e^2c \right)$. Through the isomorphism \\
 $\displaystyle\mathcal{J}: \ \ J_{Q}(X, Z)\longrightarrow E_{w}(x, y)_{\mid_{\mathbb{Q}}}$, $\left(X, \ Z \right) \longmapsto \ \left(\displaystyle \frac{2e(Z+e)+dX}{X^2}, \  \displaystyle\frac{4e^2(Z+e)+2e(dX+cX^2)-\frac{d^2X^2}{2e}}{X^3}\right)$
  and its inverse $\displaystyle \mathcal{J}^{-1}: \ \ E_{w}(x, y)_{\mid_{\mathbb{Q}}} \ \longrightarrow \ J_{Q}(X, Z)$,  $\left(x, y\right)\longmapsto \ \displaystyle\left(\frac{2e(x+c)-\frac{d^2}{2e}}{y}, \ -e+\frac{X(Xx-d)}{2e} \right)$.
}
\end{Theorem}
\begin{Proof} \ \ \\ \rm{	
$J^{-1}(x, y)=J_{Q}\displaystyle\left(X, -e+\frac{X(Xx-d)}{2e} \right)$ with $X=\displaystyle\frac{2e(x+c)-\frac{d^2}{2e}}{y}$
\begin{align*}
\displaystyle\left[-e+\frac{X(Xx-d)}{2e} \right]^2&=aX^4+bX^3+cX^2+dX+e^2 \ \ \ \ \ \ \ \ \ \ \ \ \ \ \ \ \ \ \ \ \ \ \ \ \ \ \ \ \ \ \ \ \ \ \ \ \ \ \ \ \ \ \ \ \ \ \ \ \ \ \ \ \ \ \ \ \ \ \ \ \ \ \ \ \ \ \ \ \ \ \ \ \ \ \ \ \ \ \ \ \ \ \\
\text{Since } \displaystyle\left[-e+\frac{X(Xx-d)}{2e} \right]^2 &= e^2+ \displaystyle\frac{X^2(X^2x^2-2dXx+d^2)}{4e^2}-X^2x+dX \ \text{ , then }\\
e^2+ \displaystyle\frac{X^2(X^2x^2-2dXx+d^2)}{4e^2}&-X^2x+dX = aX^4+bX^3+cX^2+dX+e^2\\
\displaystyle\frac{X^4x^2-2dX^3x+d^2X^2}{4e^2}&=aX^4+bX^3+(c+x)X^2\\
(x^2-4ae^2)X^4-(2dx+4be^2)&X^3+(d^2-4(c+x)e^2)X^2=0\\
(x^2-4ae^2)X^2-(2dx+4be^2)&X+(d^2-4(c+x)e^2)=0, \ \ \ X^2 \neq 0 \ \ \ \text{with } X=\displaystyle\frac{2e(x+c)-\frac{d^2}{2e}}{y}\\
(x^2-4ae^2)\displaystyle\left(\frac{2e(x+c)-\frac{d^2}{2e}}{y}\right) ^2&-(2dx+4be^2)\displaystyle\left(\frac{2e(x+c)-\frac{d^2}{2e}}{y}\right)+(d^2-4(c+x)e^2)=0\\
\displaystyle\frac{1}{4e^2y^2}\left(x^2-4ae^2\right)\left(4e^2(x+c)-d^2 \right)^2&-\displaystyle\frac{1}{2ey}\left(2dx+4be^2 \right)\left(4e^2(x+c)-d^2 \right)+d^2-4(c+x)e^2=0\\
\left(x^2-4ae^2\right)\left(4e^2(x+c)-d^2 \right)^2&-2ey\left(2dx+4be^2 \right)\left(4e^2(x+c)-d^2 \right)+4e^2y^2\left(d^2-4(c+x)e^2 \right)=0\\
\left(x^2-4ae^2\right)\left(4e^2(x+c)-d^2 \right)^2&-2ey\left(2dx+4be^2 \right)\left(4e^2(x+c)-d^2 \right)-4e^2y^2\left(4(c+x)e^2-d^2 \right)=0\\
\left(x^2-4ae^2\right)\left(4e^2(x+c)-d^2 \right)&-2ey\left(2dx+4be^2 \right)-4e^2y^2=0 , \ \ \ \ \left(4e^2(x+c)-d^2 \right) \neq 0 \\
4e^2x^3+(4ce^2-d^2)x^2-16ae^4x&-4ae^2(4ce^2-d^2)-2ey(2dx+4be^2)-4e^2y^2=0 \\
x^3+\displaystyle\frac{4ce^2-d^2}{4e^2}x^2-4ae^2x-a(4ce^2&-d^2)-\displaystyle\frac{y}{e}(dx+2be^2)-y^2=0 \\
x^3+\displaystyle\frac{4ce^2-d^2}{4e^2}x^2-4ae^2x-a(4ce^2&-d^2)-y(\displaystyle\frac{d}{e}x+2be)-y^2=0 \\
\text{Then} \ \ y^2+\displaystyle\frac{d}{e}xy+2bey=x^3&+\displaystyle\frac{4ce^2-d^2}{4e^2}x^2-4ae^2x-a(4ce^2-d^2)\\
\text{Hence} \ \ y^2+\displaystyle\frac{d}{e}xy+2eby=x^3&+\displaystyle\left(c- \frac{d^2}{4e^2}\right)x^2-4e^2ax+a(d^2-4ce^2)\\
& \text{By identification with respect to the Weierstrass form, we have: }  \\
a_{1}=\displaystyle\frac{d}{e}, & \ a_{3}=2eb, \ a_{2}=\displaystyle c-\frac{d^2}{4e^2}, \ a_{4}=-4e^2a, \ a_{6}=\displaystyle a\left(d^2-4e^2c \right)
\end{align*}

Now, let's prove that $\mathcal{J}^{-1}$ is effectively the inverse map of $\mathcal{J}$.
\begin{align*}
\left(\mathcal{J}\circ\mathcal{J}^{-1}\right)(x, y)&=\mathcal{J}\left(\mathcal{J}^{-1}(x, y)\right)\\
&=\mathcal{J}\left(X, \ -e+\frac{X(Xx-d)}{2e}\right)\\
&=\left(\frac{2e\left[-e+\frac{X(Xx-d)}{2e}+e \right]+dX}{X^2}, \ \frac{4e^2\left[-e+\frac{X(Xx-d)}{2e}+e \right]+2e(dX+cX^2)-\frac{d^2X^2}{2e}}{X^3}\right) \\
&=\left(\frac{X^2x-Xd+dX}{X^2}, \ \frac{2eX(Xx-d)+2edX+2ecX^2-\frac{d^2X^2}{2e}}{X^3} \right)\\
&=\left(x, \ \frac{-2eX(Xx-d)+2edX+2ecX^2-\frac{d^2X^2}{2e}}{X^3}\right) \\
&=\left(x, \ \frac{2eX^2(x+c)-\frac{d^2X^2}{2e}}{X^3} \right)\\
&=\left(x, \ \frac{4e^2X^2(x+c)-d^2X^2}{2eX^3} \right) \\
&=\left(x, \ \frac{4e^2(x+c)-d^2}{2eX} \right) \\
&=\left(x, \ \frac{4e^2y(x+c)-yd^2}{2e(2e(x+c)-\frac{d^2}{2e})} \right)\\
&=\left(x, \ y\frac{4e^2(x+c)-d^2}{4e^2(x+c)-d^2} \right) \\
&=(X, Z) \\
&\text{ This shows that } \left(\mathcal{J}\circ\mathcal{J}^{-1}\right)=id \\
&\text{On the other hand: }\\
\left(\mathcal{J}^{-1}\circ\mathcal{J}\right)(X, Z)&=\mathcal{J}^{-1} \left(\mathcal{J}(X, Z)\right) \\
&=\displaystyle\left(\frac{2e(Z+e)+dX}{X^{2}}, \ \frac{4e^2(Z+e)+2e(dX+cX^2)-\frac{d^2X^2}{2e}}{X^3} \right) \\
&=\displaystyle\left(\frac{2e\left(\frac{2e(Z+e)+dX}{X^2}+c \right)-\frac{d^2}{2e} }{\displaystyle\frac{4e^2(Z+e)+2e(dX+cX^2)-\frac{d^2X^2}{2e}}{X^3}}, \ -e+\frac{X\left(X\left(\frac{2e(Z+e)+dX}{X^2} \right)-d  \right) }{2e}\right)\\
&=\displaystyle\left(\frac{2e\left(\frac{2e(Z+e)+dX+cX^2}{X^2} \right)-\frac{d^2}{2e} }{\displaystyle\frac{4e^2(Z+e)+2e(dX+cX^2)}{X^3}-\frac{d^2X^2}{2eX^3}}, \ -e+\frac{2e(Z+e)+dX-dX}{2e } \right) \\
&=\displaystyle\left(\frac{\displaystyle\frac{4e^2(Z+e)+2edX+2ecX^2}{X^2}-\frac{d^2}{2e}}{\displaystyle\frac{4e^2(Z+e)+2e(dX+cX^2)}{X^3}-\frac{d^2X^2}{2eX^3}}, \ -e+\frac{2e(Z+e)}{2e} \right) \\
&=\displaystyle\left(\displaystyle\frac{\displaystyle\frac{4e^2(Z+e)+2e(dX+cX^2)}{X^2}-\frac{d^2}{2e}}{\displaystyle\frac{4e^2(Z+e)+2e(dX+cX^2)}{X^3}-\frac{d^2}{2eX}}, \ -e+Z+e \right) \\
&=\displaystyle\left(\frac{\displaystyle\frac{8e^3(Z+e)+4e^2(dX+cX^2)-d^2X^2}{2eX^2}}{\displaystyle\frac{8e^3X(Z+e)+4e^2X(dX+cX^2)-d^2X^3}{2eX^4}}, \ \ Z \right) \\
&=\displaystyle\left(\frac{\displaystyle\frac{8e^3(Z+e)+4e^2(dX+cX^2)-dX^2}{2eX^2}}{\displaystyle\frac{8e^3(Z+e)+4e^2(dX+cX^2)-d^2X^2}{2eX^3}}, \ \ Z \right)  \\
&=\displaystyle\left(\frac{8e^3(Z+e)+4e^2(dX+cX^2)-d^2X^2}{\frac{8e^3(Z+e)+4e^2(dX+cX^2)-d^2X^2}{X}}, \ \ Z \right)\\
&=(x, \ \ y)\\
&\text{ This also shows that } \left(\mathcal{J}^{-1}\circ\mathcal{J}\right)=id, \\
& \text{Hence $J$ is a bijective morphism since $\left(\mathcal{J}^{-1}\circ\mathcal{J}\right)=\left(\mathcal{J}\circ\mathcal{J}^{-1}\right)=id$, thus and isomorphism}
\end{align*}
}
\end{Proof}

\begin{Proposition} \ \\ \label{prop27} \rm{
	The equation $E_{w}(x, y)_{\mid_{\mathbb{Q}}}$ is isomorphic to the short Weierstrass form curve $E_{k \neq 2, 3}(X^{'}, Y): \ Y^2=X^{'3}+aX^{'}+b$ in characteristic $\neq 2, 3$ through the isomorphism 
	$\mathcal{K}: \ E_{w}(x, y)_{\mid_{\mathbb{Q}}} \longrightarrow E_{k \neq 2, 3}(X^{'}, Y)$, $(x, y) \longmapsto \ (x+\frac{A}{3}, \ y+\frac{a_{1}x+a_{3}}{2})$
	where $a=B-\frac{A^{2}}{3}$ and $b=\frac{27C-9AB+2A^{3}}{27}$ with $A=a_{2}+\frac{a_{1}^{2}}{4}$, $B=a_{4}+\frac{a_{1}a_{3}}{2}$, $C=a_{6}+\frac{a_{3}^{2}}{4}$, $Y=y+\frac{a_{1}x+a_{3}}{2}$, $X^{'}=x+\frac{A}{3}$. \\
	Here $\mathcal{K}^{-1}: E_{k \neq 2, 3}(X^{'}, Y) \longrightarrow  E_{w}(x, y)_{\mid_{\mathbb{Q}}}$, $(X^{'}, Y) \longmapsto \ (X^{'}-\frac{A}{3}, \ Y-\frac{a_{1}x+a_{3}}{2})$.
}
\end{Proposition}
\begin{Proof} \ \\ \rm{
Considering
$E_{w}(x, y): \ y^2+a_{1}xy+a_{3}y=x^3+a_{2}x^2+a_{4}x+a_{6}$.\\
$\left(y+\frac{a_{1}x+a_{3}}{2}\right)^2=x^3+\left(a_{2}+\frac{a_{1}^{2}}{4}\right)x^2+\left(a_{4}+\frac{a_{1}a_{3}}{2}\right)x+a_{6}+\frac{a_{3}^{2}}{4}$. Setting 
$Y=y+\frac{a_{1}x+a_{3}}{2}$ with 
$A=a_{2}+\frac{a_{1}^{2}}{4}$, $B=a_{4}+\frac{a_{1}a_{3}}{2}$ and $C=a_{6}+\frac{a_{3}^{2}}{4}$.\\
$Y^2=(x+\frac{A}{3})^3+(B-\frac{A^2}{3})(x+\frac{A}{3})+\frac{27C-9AB+2A^{3}}{27}$. Setting $X^{'}=x+\frac{A}{3}$,\\
 the equation becomes $E_{k \neq 2, 3}(X^{'}, Y): \ Y^2=X^{'3}+aX^{'}+b$, With $a=B-\frac{A^2}{3}$ and $b=\frac{27C-9AB+2A^{3}}{27}$.\\ 
 $\displaystyle\left( \mathcal{K}\circ\mathcal{K}^{-1}\right)=\displaystyle\left( \mathcal{K}^{-1}\circ\mathcal{K}\right)=\ id$ is straightforward.
}
\end{Proof}
}
\begin{Theorem} \ \label{ntsukunu} \\ \rm{
$\forall \ (X^{'}, Y)\in E_{k \neq 2, 3}(X^{'}, Y)$, set $x=X^{'}-\frac{A}{3}$ and $y=Y-\frac{a_{1}x+a_{3}}{2}$. \\
$(X, \ Z)=\left(\frac{2e(x+c)-\frac{d^2}{2e}}{y}, \ -e+\frac{X(Xx-d)}{2e} \right)$.\\
Setting $Y=\displaystyle\frac{1}{[32(n^3-n)^2]}\left( 32n^3(n^2+1)Z- \left[64n^{4}X^{2}-(384n^5+64n^7+64n^3)X+(192n^8+640n^6+192n^4) \right]\right)$, \\
If $\left(X, \ Y \right)\in \mathbb{Z}_{>0}\times \mathbb{Z}_{>0}$  then $(X, \ Y) \in R_{n}(X, Y)$	and a prime factorization of $n$ is found.
	
}
\end{Theorem}
\begin{Proof}  \ \\ \rm{
We give here a proof in 4 backward steps.
Let's consider Proposition \ref{prop27}\\
Consider $\mathcal{K}: \ E_{w}(x, y)_{\mid_{\mathbb{Q}}} \longrightarrow \ E_{k \neq 2, 3}(X^{'}, Y)$, $(x, \ y) \longmapsto \ \left(x+\frac{A}{3}, \ y+\frac{a_{1}x+a_{3}}{2} \right)$. \\
First backward step: $\forall \ (X^{'}, Y)\in E_{k \neq 2, 3}(X^{'}, Y)$, $\mathcal{K}^{-1}(X^{'}, Y)=(x, \ y)\in E_{w}(x, y)_{\mid_{\mathbb{Q}}}$. \\
Consider $\mathcal{J}: \ J_{Q}(X, Z)\longrightarrow E_{w}(x, y)_{\mid_{\mathbb{Q}}}$, $\left(X, \ Z \right) \longmapsto (x, \ y)$. \\
Second backward step: $\forall \ (x, \ y)\in E_{w}(x, y)_{\mid_{\mathbb{Q}}}$, $\mathcal{J}^{-1}(x, \ y)=\mathcal{J}^{-1}\left(\mathcal{K}^{-1}(X^{'}, Y)\right)=\left(\mathcal{J}^{-1}\circ\mathcal{K}^{-1}\right)(X^{'}, Y)=(X, \ Z)\in J_{Q}(X, Z)$.\\
Consider $\mathcal{I}: \  R_{n}(X, Y) \longrightarrow \ J_{Q}(X, Z)$, $\left(X, Y \right) \longmapsto (X, Z)$.\\
Third backward step: $\mathcal{I}^{-1}(X, Z)=\mathcal{I}^{-1}\left( \mathcal{J}^{-1}(\mathcal{K}^{-1}(X^{'}, Y))\right)=\left(\mathcal{I}^{-1}\circ \mathcal{J}^{-1}\circ\mathcal{K}^{-1}\right)(X^{'}, Y)=(X, \ Y)\in R_{n}(X, Y)$.\\
Since $Z=\frac{(n^2-1)^2}{n(n^2+1)}\left( Y+\frac{64n^{4}X^{2}-(384n^{5}+64n^{7}+64n^{3})X+(192n^{8}+640n^{6}+192n^{4})}{32(n^3-n)^2}\right)$, then \\
$Y=\displaystyle\frac{1}{[32(n^3-n)^2]}\left( 32n^3(n^2+1)Z- \left[64n^{4}X^{2}-(384n^5+64n^7+64n^3)X+(192n^8+640n^6+192n^4) \right]\right)$. If $Y\in \mathbb{Z}_{+}$, then $(X, Y)\in R_{n}(X, Y)$.\\
Consider: $h: \lbrace P_{2}, P_{3}\rbrace \subset \mathcal{B}_{n}(x, y)_{\mid_{x\geq 4n}} \longrightarrow \ R_{n}(X, Y)$, $\lbrace P_{2}, P_{3}\rbrace\longmapsto \  \left(x_{P_{2}}+x_{P_{3}}, \ x_{P_{2}}x_{P_{3}}\right)$.\\
Fourth backward step: $h^{-1}\left(X, Y\right)=h^{-1}\left(\mathcal{I}^{-1}\left( \mathcal{J}^{-1}(\mathcal{K}^{-1}(X^{'}, Y))\right)\right)=\left(h^{-1}\circ\mathcal{I}^{-1}\circ \mathcal{J}^{-1}\circ\mathcal{K}^{-1}\right)(X^{'}, Y)=\lbrace P_{2}, P_{3} \rbrace \rbrace \subset \mathcal{B}_{n}(x, y)_{\mid_{x\geq 4n}}$ where $P_{2}=(x_{P_{2}}, \sqrt{x_{P_{2}}^2-4nx_{P_{2}}})$ and $P_{3}=(x_{P_{3}}, \sqrt{x_{P_{3}}^2-4nx_{P_{3}}})$. Hence, a prime factorization of $n$ is found.
}
\end{Proof}
\begin{Remark} \ \\ \rm{
This proposed factorization approach is summarized by the following diagram: \\ \\
\begin{xy}
	(0,20)*+{\lbrace P_{2}, P_{3}\rbrace \subset \mathcal{B}_{n}(x, y)_{\mid_{x\geq 4n}}}="a"; (40,20)*+{ \ R_{n}(X, Y)}="b"; (70,20)*+{\ \ \ \ J_{Q}(X, Z)}="c";(102,20)*+{ \ \ \ \ \ \ E_{w}(x, y)_{\mid_{\mathbb{Q}}}}="d"; (142,20)*+{ \ \ \ \ \ \ \ \ \ \ \ E_{k \neq 2, 3}(X^{'},Y) \ \ }="e"; 
	{\ar@{->}@/^{4pc}/ "a";"e"}?*!/_4mm/{ \mathcal{K}\circ \mathcal{J} \circ \mathcal{I} \circ h};
	{\ar@<1.ex>"a";"b"}?*!/_2mm/{\ \ \ \ \ \ \ \ \ h};%
	{\ar@<1.ex> "b";"a"}?*!/_2mm/{\ \ \ \ \ \ \ \ \ h^{-1}};
	{ \ \ \ar@<1.ex>"b";"c"}?*!/_2mm/{\ \mathcal{I}};%
	{\ar@<1.ex> "c";"b"}?*!/_2mm/{\ \mathcal{I}^{-1}};
	{\ar@<1.ex>"c";"d"}?*!/_2mm/{\ \mathcal{J}};%
	{\ar@<1.ex> "d";"c"}?*!/_2mm/{\ \mathcal{J}^{-1}};
	{\ar@<1.ex>"d";"e"}?*!/_2mm/{\mathcal{K} \ \ \ \ };%
	{\ar@<1.ex>"e";"d"}?*!/_2mm/{\mathcal{K}^{-1} \ \ \ };
	{\ar@{->}@/^{4pc}/ "e";"a"}?*!/_5mm/{ h^{-1}\circ\mathcal{I}^{-1}\circ \mathcal{J}^{-1} \circ \mathcal{K}^{-1}}; %
\end{xy}
\ \ \\ \\
Under the condition of Theorem \ref{ntsukunu}, 
$\forall \ (X^{'}, Y) \in  E_{w}(x, y)_{\mid_{\mathbb{Q}}}$, \ $\left( h^{-1}\circ\mathcal{I}^{-1}\circ \mathcal{J}^{-1} \circ \mathcal{K}^{-1}\right) (X^{'}, Y)=\left(x_{P_{2}}+x_{P_{3}}, \ x_{P_{2}}x_{P_{3}} \right)$ where $\lbrace P_{2}, P_{3}\rbrace \subset \mathcal{B}_{n}(x, y)_{\mid_{x\geq 4n}}$.

}	
\end{Remark}
\textbf{Comment on $E_{w}(x, y)_{\mid_{\mathbb{Q}}}$ rational torsion points} \\ \rm{
We first note that  $E_{w}(x, y)_{\mid_{\mathbb{Q}}}$ has a 2-torsion rational point $P=(X^{'}, Y)=(\frac{1}{3}(-2a_{2}+\frac{a_{1}^{2}}{4}), \frac{1}{2}(a_{1}a_{2}-a_{3}))$.\\
To prove this, let's first consider,  $E_{k \neq 2, 3}(X^{'}, Y): \ Y^2=X^{'3}+aX^{'}+b$, With $a=B-\frac{A^2}{3}$ and $b=\frac{27a_{6}-9AB+2A^{3}}{27}$.\\
From the Jacobi quartic, setting $X=0$, we get the two points $(X, Z_{1})=(0, e)$ and $(X, Z_{2})=(0, -e)$\\
And through the isomorsphism $J$, we have: \\
$\mathcal{J}(X, Z_{1})=\mathcal{J}(0, e)=(0:1:0)=\infty$ and $\mathcal{J}(X, Z_{2})=\mathcal{J}(0, -e)$ \\
 $X^{'}=x+\frac{A}{3}=-a_{2}+\frac{1}{3}(a_{2}+\frac{a_{1}^{2}}{4})=\frac{1}{3}(-2a_{2}+\frac{a_{1}^{2}}{4})$ and $Y=y+\frac{a_{1}x+a_{3}}{2}=a_{1}a_{2}-a_{3}+\frac{-a_{1}a_{2}+a_{3}}{2}=\frac{1}{2}(a_{1}a_{2}-a_{3})$.
hence a rational point is $(X^{'}, Y)=(\frac{1}{3}(-2a_{2}+\frac{a_{1}^{2}}{4}), \frac{1}{2}(a_{1}a_{2}-a_{3}))$. \\
}

\section*{\textbf{Conclusion}}  \rm{
In this paper, we present an approach to factoring integers of type of RSA modulus using arithmetical properties of a hyperbola. This approach aims at finding points $P_{2}$ and $P_{3}$ of $ \mathcal{B}_{n}(X, Y, Z)_{\mid_{x\geq 4n}}= \displaystyle \lbrace \left(X: Y: Z\right)\in \mathbb{P}^{2}(\mathbb{Q}) \ / \ \displaystyle Y^{2}=X^{2}-4nXZ \rbrace$ whose knowledge leads the factorization of the RSA modulus $n$. We prove that this solution can be found on a singular weierstrass curve that is isomorphic to a Jacobi quartic with the introduction of both concepts of Hyperbola X-root and Hyperbola Y-root.\\ 
As perspective, we shall investigate its performance and complexity over rationals and finite fields. We should also generate a dataset of enough $n$, $\varepsilon$ with related features to apply machine learning models for pattern investigation and prediction. 
}

\end{document}